\documentclass{amsart}

\theoremstyle{plain}
\newtheorem{theorem}{Theorem}[section]
\newtheorem{proposition}[theorem]{Proposition}

\newtheorem{corollary}[theorem]{Corollary}
\newtheorem{conjecture}[theorem]{Conjecture}

\theoremstyle{definition}

\newtheorem{example}[theorem]{Example}
\newtheorem{remark}[theorem]{Remark}

\begin{document}

\title{A note on a degree sum condition for long cycles in graphs}
\author{Janusz Adamus}
\address{Department of Mathematics, Middlesex College, University of Western Ontario, London, Ontario N6A 5B7 Canada}
\email{jadamus@uwo.ca}

\begin{abstract}
We conjecture that a $2$-connected graph $G$ of order $n$, in which $d(x)+d(y)\geq n-k$ for every pair of non-adjacent vertices $x$ and $y$, contains a cycle of length $n-k$ \ ($k<n/2$), unless $G$ is bipartite and  $n-k$ is odd. This generalizes to long cycles a well-known degree sum condition for hamiltonicity of Ore. The conjecture is shown to hold for $k=1$.
\end{abstract}
\maketitle

\section{Introduction}
\label{sec:intro}

The subject of this note is the following conjecture, in which we generalize to long cycles a well-known degree sum condition for hamiltonicity of Ore~\cite{Ore}. All graphs considered are finite, undirected, with no loops or multiple edges.
\smallskip

\begin{conjecture}
\label{conj:long-cycles}
Let $G$ be a $2$-connected graph of order $n\geq3$, $n\neq5,7$, and let $k<n/2$ be an integer.
If
\[
d(x)+d(y)\geq n-k
\]
for every pair of non-adjacent vertices $x$ and $y$, then $G$ contains a cycle of length $n-k$, unless $G$ is bipartite and \ $n-k\equiv1\ (\!\!\!\!\mod2)$.
\end{conjecture}
\medskip

\begin{remark}
\label{rem:sharp}
The conjecture is sharp. First of all, a quick look at $C_5$ and $C_7$ ensures that the assumption $|G|\neq5,7$ is necessary. Secondly, it is easy to see that without the $2$-connectedness assumption, there could be no long cycles at all. Consider, for instance, a graph $G$ obtained from disjoint cliques $H_1=K_{\left\lfloor n/2\right\rfloor}$ and $H_2=K_{\left\lceil n/2\right\rceil}$ by joining a single vertex $x_0$ of $H_2$ with every vertex of $H_1$.
Finally, the bound for the degree sum of non-adjacent vertices is best possible, as shown in the example below.
\end{remark}
\medskip

\begin{example}
\label{ex:sharp}
Let $G$ be a graph obtained from the complete bipartite graph $K_{(n-k-1)/2,(n+k+1)/2}$ by joining all the vertices in the smaller colour class. Then $d(x)+d(y)\geq n-k-1$ for every pair of non-adjacent vertices $x$ and $y$, and $G$ contains no cycle of length greater than $n-k-1$.
\end{example}
\medskip

Our main result is the following theorem that implies Conjecture~\ref{conj:long-cycles} for $k=1$, as shown in Section~\ref{sec:long-cycles}. The proof of Theorem~\ref{thm:main} is given in the last section.
\smallskip

\begin{theorem}
\label{thm:main}
Let $G$ be a $2$-connected graph of order $n\geq3$, in which 
\[
d(x)+d(y)\geq n-1
\]
 for every pair of non-adjacent vertices $x$ and $y$.
\begin{itemize}
\item[(i)] If $n$ is even, then $G$ is hamiltonian.
\item[(ii)] If $n$ is odd, then $G$ contains a cycle of length at least $n-1$. 
\end{itemize}
Moreover, $G$ is not hamiltonian only if the minimal degree of its $n$-closure, $Cl_n(G)$, equals $(n-1)/2$. In this case, $Cl_n(G)$ is a maximal non-hamiltonian graph.
\end{theorem}
\medskip

Recall that the \emph{$n$-closure} $Cl_n(G)$ of $G$ is a graph obtained from $G$ by succesively joining all pairs $(x,y)$ of non-adjacent vertices satisfying $d(x)+d(y)\geq n$.

\medskip

\section{Long cycles in graphs}
\label{sec:long-cycles}

\begin{proposition}
\label{prop:conj-holds}
Conjecture~\ref{conj:long-cycles} holds for $k=1$.
\end{proposition}
\medskip

For the proof, we will need the following result of~\cite{H-F-S}:
\smallskip

\begin{theorem}[Haggkvist-Faudree-Schelp]
\label{thm:HFS}
Let $G$ be a hamiltonian graph on $n$ vertices. If $G$ contains more than $\left\lfloor\frac{(n-1)^2}{4}\right\rfloor+1$ edges, then $G$ is pancyclic or bipartite.
\end{theorem}
\bigskip

\textit{Proof of Proposition~\ref{prop:conj-holds}.}
By Theorem~\ref{thm:main}, we may assume that $G$ is hamiltonian. Suppose first that $G$ is a $2$-connected non-bipartite hamiltonian graph of order $n$, in which $d(x)+d(y)\geq n-1$ whenever $xy\notin E(G)$.\par

Consider a vertex $x$ of minimal degree $d(x)=\delta(G)$ in $G$. Write $\delta=\delta(G)$. Then $G$ has precisely $n-1-\delta$ vertices non-adjacent to $x$, each of degree at least $n-1-\delta$. The remaining $\delta+1$ vertices are of degree at least $\delta$ each, hence
\[
\left\|G\right\|\geq\frac{1}{2}[(\delta+1)\delta+(n-1-\delta)^2]\,.
\]
As $\delta\geq2$, one immediately verifies that
\[
\frac{1}{2}[(\delta+1)\delta+(n-1-\delta)^2]>\frac{(n-1)^2}{4}+1\,,
\]
whenever $n\neq5$.\par

It remains to consider the case of $G$ a bipartite $2$-connected hamiltonian graph of order $n$. But then $n$ must be even, for otherwise $G$ would contain an odd cycle. Thus $n-1\equiv1\ (\!\!\!\mod2)$, which completes the proof.\qed
\medskip

For convenience, let us finally recall two well-known results, that we shall need in the proof of Theorem~\ref{thm:main}:
\smallskip

\begin{theorem}[Dirac~\cite{Dirac}]
\label{thm:I}
Let $G$ be a graph of order $n\geq3$ and minimal degree $\delta(G)\geq n/2$. Then $G$ is hamiltonian.
\end{theorem}
\smallskip

\begin{theorem}[Bondy-Chvatal~\cite{B-C}]
\label{thm:III}
Let $G$ be a graph of order $n$ and suppose that there is a pair of non-adjacent vertices $x$ and $y$ of $G$ such that $d(x)+d(y)\geq n$. Then $G$ is hamiltonian if and only if $G+xy$ is hamiltonian.
\end{theorem}
\smallskip

\begin{corollary}
\label{thm:V}
A graph $G$ is hamiltonian if and only if its $n$-closure $Cl_n(G)$ is so.
\end{corollary}

\medskip

\section{Proof of Theorem~\ref{thm:main}}
\label{sec:proof-main-theorem}

\subsubsection*{Proof of part (i)} Suppose there exists an even integer $n\geq4$ for which the assertion of the theorem does not hold. Let $G$ be a maximal non-hamiltonian $2$-connected graph of order $n$, in which $d(x)+d(y)\geq n-1$ whenever $xy\notin E(G)$.\par

By maximality of $G$, $G+xy$ is hamiltonian for every pair of non-adjacent vertices $x,y\in V(G)$. Hence, by Theorem~\ref{thm:III}, we must have
\[
\tag{$\ast$}
d(x)+d(y)=n-1 \quad\textrm{whenever}\ xy\notin E(G)\,.
\]
The minimal degree $\delta(G)$ of $G$ satisfies inequality $\delta(G)<n/2$, by Theorem~\ref{thm:I}, hence, in particular, $n-1-\delta(G)\geq\delta(G)+1$.
\smallskip

Pick $x\in V(G)$ with $d(x)=\delta(G)$. There are precisely $n-1-\delta(G)$ vertices in $G$ non-adjacent to $x$, each of degree $n-1-\delta(G)$, by $(\ast)$. Put $V=\{v\in V(G):xv\notin E(G)\}$. Pick $y\in V$. As $d(y)=n-1-\delta(G)$, there are precisely $\delta(G)$ vertices in $G$ non-adjacent to $y$, each of degree $\delta(G)$, by $(\ast)$ again. Put $U=\{u\in V(G):uy\notin E(G)\}$. Then $|U|=\delta(G)$, $|V|=n-1-\delta(G)$, and $U\cap V=\emptyset$, because vertices in $U$ are of degree $\delta(G)$ and those in $V$ are of degree $n-1-\delta(G)>\delta(G)$. It follows that there exists a vertex $z$ in $G$ such that $V(G)=U\cup V\cup \{z\}$ is a partition of the vertex set of $G$.
\smallskip

We will now show that $d(z)=n-1$:
Observe first that $d(z)>\delta(G)$. Indeed, if $d(z)=\delta(G)$, then by $(\ast)$, $z$ is adjacent to every vertex in $U$, as $2\delta(G)<n-1$. But $z$ is also adjacent to $y$, as $z\notin U$, hence $d(z)\geq|U|+1=\delta(G)+1$; a contradiction. Consequently, $z$ is adjacent to every vertex in $V$, by $(\ast)$ again, as $d(z)+(n-1-\delta(G))>n-1$. Hence $d(z)\geq|V|=n-1-\delta(G)$. On the other hand, $z$ is adjacent to $x$, as $z\notin V$, which yields $d(z)\geq|V|+1=n-\delta(G)$. This last inequality paired with $(\ast)$ implies that $z$ is adjacent to every other vertex in $G$, as required.
\smallskip

Next observe that $u_1u_2\in E(G)$ for every pair of vertices $u_1,u_2$ in $U$, as $d(u_1)+d(u_2)=2\delta(G)<n-1$. It follows that $N(u)\supset U\cup\{z\}\setminus\{u\}$, and hence, by comparing cardinalities, $N(u)=U\cup\{z\}\setminus\{u\}$ for every $u\in U$.\par

Similarly, $v_1v_2\in E(G)$ for every pair $v_1,v_2$ in $V$, hence $N(v)=V\cup\{z\}\setminus\{v\}$ for every $v\in V$. Therefore $G=G_1\cup G_2$, where $G_1$ is a complete graph of order $\delta(G)+1$ spanned on the vertices of $U\cup\{z\}$, and $G_2$ is a complete graph of order $n-\delta(G)$ spanned on $V\cup\{z\}$. Then $z$ is a cutvertex, contradicting the assumption that $G$ be $2$-connected.
\medskip

\subsubsection*{Proof of part (ii)} Suppose there exists a $2$-connected graph of odd order $n\geq3$, in which $d(x)+d(y)\geq n-1$ for every pair of non-adjacent vertices $x$ and $y$, that does not contain neither a Hamilton cycle nor a cycle of length $n-1$. Let $G$ be maximal such a graph of order $n$. By maximality of $G$, $G+xy$ contains a cycle of length at least $n-1$ whenever $xy\notin E(G)$. Hence $G$ contains a path of length at least $n-2$ between any two of its non-adjacent vertices.\par

Pick a pair of non-adjacent vertices $x$ and $y$. By a theorem of P{\'o}sa, $G$ contains a Hamilton $x-y$ path $P$, and hence, by Theorem~\ref{thm:III}, the sum $d(x)+d(y)$ actually equals $n-1$. Write $P=u_1u_2\dots u_n$, where $u_1=x$ and $u_n=y$.\par

Put $I_x=\{i:xu_{i+1}\in E(G),1\leq i\leq n-1\}$ and $I_y=\{i:u_iy\in E(G),1\leq i\leq n-1\}$.
If $I_x\cap I_y\neq\emptyset$, say $i_0\in I_x\cap I_y$, then $G$ contains a Hamilton cycle
\[
u_1u_{i_0+1}u_{i_0+2}\dots u_nu_{i_0}u_{i_0-1}\dots u_2u_1\,.
\]
We may thus assume that $I_x\cap I_y=\emptyset$. Then, for every $1\leq i\leq n-1$, either $u_i$ is adjacent to $y$ or else $u_{i+1}$ is adjacent to $x$, because $|I_x|+|I_y|=d(x)+d(y)=n-1$. Let $d=d(y)$ and let $v_1,\dots,v_d\!=\!y$ be the vertices that lie on $P$ next to the (respective) neighbours of $y$.\par

If there exists $j<d$ such that $v_j\notin N(y)$, then $v_j=u_{i_0}$ for some $i_0\in I_x$. It follows that $u_{i_0+1}$ is adjacent to $x$, and $G$ contains a cycle of length $n-1$ of the form
\[
u_1u_{i_0+1}u_{i_0+2}\dots u_nu_{i_0-1}u_{i_0-2}\dots u_2u_1\,.
\]
Therefore we can assume that
\[
\tag{$\dagger$}
v_1,\dots,v_{d-1}\ \mathrm{are\ all\ adjacent\ to}\ y\,.
\]
Let $z$ denote the furthermost neighbour of $y$ on $P$. It follows from $(\dagger)$ that all the vertices between $z$ and $y$ on $P$ are adjacent to $y$, and hence $z=u_{n-d}$.\par

Suppose $N(v_j)\subset\{z,v_1,\dots,v_d\}$ for $j\leq d$. Then $N(u_i)\subset\{u_1,\dots,u_{n-d-1}, z\}$ for $i\leq n-d-1$. Consequently, $d(u_i)\leq n-d-1$, $d(v_j)\leq d$, and $u_iv_j\notin E(G)$ for $i\leq n-d-1$ and $j\leq d$. But then $d(u_i)+d(v_j)\geq n-1$ yields
\[
d(u_i)=n-d-1 \quad\mathrm{and}\quad d(v_j)=d\quad\mathrm{for\ \ }i=1,\dots,n-d-1,\ j=1,\dots,d\,.
\]
Therefore, as in the proof of part \emph{(i)}, we get that $G=G_1\cup G_2$, where $G_1$ is a complete graph of order $n-d$ spanned on the vertices $\{u_1,\dots,u_{n-d-1},z\}$ and $G_2$ is a complete graph of order $d+1$ on $\{z,v_1,\dots,v_d\}$. Then $z$ is a cutvertex contradicting our assumptions on $G$.
\smallskip

It remains to consider the case of some $v_{j_0}$ being adjacent to $u_{i_0}$, where $i_0\leq n-d-1$. But then again $G$ contains a Hamilton cycle
\[
u_1\dots u_{i_0}v_{j_0}\dots v_dv_{j_0-1}\dots u_{i_0+1}u_1\,.
\]
\bigskip

For the proof of the last assertion of Theorem~\ref{thm:main}, suppose that $n=2k+1$ is odd and $G$ is a non-hamiltonian $2$-connected graph on $n$ vertices, satisfying $d(x)+d(y)\geq n-1$ for every pair of non-adjacent $x$ and $y$. Then the $n$-closure of $G$, $G^*=Cl_n(G)$ is not hamiltonian either, by Theorem~\ref{thm:V}, and we have equality
\[
d_{G^*}(x)+d_{G^*}(y)=n-1 \quad\textrm{whenever}\quad xy\notin E(G^*). 
\]
Now, if $\delta(G^*)<k=\frac{n-1}{2}$, then $n-1-\delta(G^*)>\delta(G^*)$ and one can repeat the proof of part \emph{(i)} to show that $G^*$ contains a Hamilton cycle, which contradicts the assumptions on $G$.\par

Thus $\delta(G^*)=\frac{n-1}{2}$. Moreover, $d_{G^*}(x)+d_{G^*}(y)=n-1=2k$ for $xy\notin E(G^*)$ implies that $d_{G^*}(x)=k$ or $d_{G^*}(x)=n-1$ for every vertex $x$.\par

Suppose $G^*$ is not maximal among the non-hamiltonian $2$-connected graphs on $n$ vertices. Then $G^*$ has a pair of non-adjacent vertices $x$ and $y$ such that $G^*+xy$ is contained in a maximal non-hamiltonian graph $H$. By maximality of $H$, $H+uv$ contains a Hamilton cycle for every $uv\notin E(H)$, so Theorem~\ref{thm:III} implies that $d_H(u)+d_H(v)=n-1$ for every $uv\notin E(H)$.\par

Notice that $d_{G^*}(x)=k$, as $d_{G^*}(x)<n-1$. Then $d_H(x)\geq k+1$ and hence, for every $v$ non-adjacent to $x$ in $G^*$, $d_H(x)+d_H(v)\geq d_{G^*}(x)+1+d_{G^*}(v)>n-1$, implying $xv\in E(H)$. Therefore $H$ is obtained from $G$ by increasing degrees of at least $x$ and all its non-neighbours in $G^*$, that is, at least $1+(n-1-k)=k+1$ vertices. But then $H$ contains at least $k+1$ vertices of degree $n-1$, which means that $\delta(H)\geq k+1=\frac{n+1}{2}$, and hence $H$ is hamiltonian by Theorem~\ref{thm:I}; a contradiction. \qed

\medskip
\bibliographystyle{amsplain}

\end{document}